\documentclass[12pt]{amsart}
\usepackage{amssymb}
\usepackage{amsbsy}
\usepackage{amscd}
\advance\textheight by \topskip
\makeatletter
\def\cal{\mathcal}
\def\Bbb{\mathbb}
\def\frak{\mathfrak}

\newenvironment{pf*}[1]{\proof[#1]}{\endproof}
\newcommand{\rom}{\textup}
\renewcommand{\thesubsection}{\thesection(\@roman\c@subsection)}
\makeatother
\newtheorem{Theorem}[equation]{Theorem}

\newtheorem{Lemma}[equation]{Lemma}
\newtheorem{Proposition}[equation]{Proposition}

\theoremstyle{definition}

\newtheorem{Example}[equation]{Example}

\makeatletter
\renewcommand\section{\@startsection{section}{1}%
  {\z@}{.7\linespacing\@plus\linespacing}{.5\linespacing}%
  {\reset@font\normalfont\bfseries\centering}}
\makeatother

\theoremstyle{remark}
\newtheorem{Remark}[equation]{Remark}

\numberwithin{equation}{section}
\numberwithin{figure}{section}

\newcommand{\thmref}[1]{Theorem~\ref{#1}}
\newcommand{\secref}[1]{\S\ref{#1}}

\newcommand{\propref}[1]{Proposition~\ref{#1}}
\newcommand{\remref}[1]{Remark~\ref{#1}}

\newcommand{\subsecref}[1]{\S\ref{#1}}

\newcommand{\Romnum}[1]{\expandafter\uppercase\expandafter{\romannumeral #1}} 
\newcommand{\C}{{\Bbb C}}
\newcommand{\Z}{{\Bbb Z}}

\newcommand{\R}{{\Bbb R}}
\newcommand{\F}{{\Bbb F}}

\newcommand{\CP}{\operatorname{\C P}}

\newcommand{\Spin}{\operatorname{\rm Spin}}
\newcommand{\Spinc}{\Spin^{c}}

\newcommand{\Map}{\operatorname{Map}}

\newcommand{\im}{\mathop{\text{\rm im}}\nolimits}

\newcommand{\ind}{\mathop{\text{\rm ind}}\nolimits}

\newcommand{\M}{{\cal M}} 
\newcommand{\G}{{\cal G}} 
\newcommand{\CC}{{\cal C}} 
\newcommand{\m}{\frak m}
\newcommand{\Diff}{\operatorname{Diff}} %

\newcommand{\SW}{\operatorname{SW}}
\newcommand{\Zp}{\Z_p}
\begin{document}
\title[Mod $p$ equality theorem]{Mod $p$ equality theorem for Seiberg-Witten invariants under $\Z_p$-actions}
\author{Nobuhiro Nakamura}
\address{Graduate School of Mathematical Sciences, University of Tokyo, 3-8-1, Komaba, Meguro-ku, Tokyo, 153-8914, Japan}
\email{nobuhiro@ms.u-tokyo.ac.jp}
\begin{abstract}
When a cyclic group $G$ of prime order acts on a $4$-manifold $X$, we prove a formula which relates the Seiberg-Witten invariants of $X$ to those of $X/G$. 
\end{abstract}
%
%
%
\maketitle
%
\section{Introduction}\label{sec:intro}
%
%

The Seiberg-Witten invariants under group actions are investigated by many authors.
In several cases, one can relate the Seiberg-Witten invariants of a $4$-manifold $X$ with an action of a group $G$ to those of its quotient ($V$-)manifold $X/G$.
In fact, in the case of free actions of prime order cyclic groups $G=\Z_p$, it is proved that the Seiberg-Witten invariant of $X$ is equal modulo $p$ to a sum of invariants of $X/G$, by Ruan-Wang~\cite{RW}, Szymik~\cite{Sz} and the author~\cite{N}.
This {\it mod $p$ equality theorem} is extended to the case of double branched coverings by Ruan-Wang~\cite{RW}, B.~D.~Park~\cite{Park} and Cho-Hong~\cite{Cho}.
On the other hand, F.~Fang~\cite{Fang} proved a {\it mod $p$ vanishing theorem} for $\Z_p$-actions. 
This is extended by the author~\cite{Np}, and in the view point there, the mod $p$ vanishing theorem can be considered as a version of mod $p$ equality theorem: If all the involved invariants of $X/G$  are $0$ by reason of negative dimensional moduli, then the invariant of $X$ is divisible by $p$.

In this paper, we shall prove a mod $p$ equality theorem for $\Z_p$-actions in  somewhat general cases. 
First, let us fix the notation. 
For an oriented closed $4$-manifold $X$ with a $\Spinc$-structure $c$, the Seiberg-Witten invariant of $(X,c)$ is denoted by $\SW(X,c)$, and the virtual dimension of the moduli by $d(c)$.
Suppose it is given an orientation-preserving action of a finite group $G$ on $X$, and the $G$-action has a lift to $c$. 
In general, there are several ways of such liftings, and we use the suffix $\alpha$ to parameterize these lifts as $G_\alpha$.
When the data $(X,c,G_\alpha)$ of a $4$-manifold $X$ with a $G$-action, a $\Spinc$-structure $c$, and a lift $G_\alpha$ of the $G$-action to $c$ are given, Y.~Ruan~\cite{Ruan} defined the $G$-monopole invariant, denoted by $\SW(X,c,G_\alpha)$, which is naturally identified with the Seiberg-Witten invariant of $V$-manifold $X/G$ with a $V$-$\Spinc$-structure $c/G_\alpha$. 
(See \secref{sec:G-monopole}.)
The virtual dimension of the moduli of $G_\alpha$-invariant solutions is denoted by $d(c,G_\alpha)$.
For any $G$-space $Y$, let $Y^G$ be the fixed point set of the $G$-action. Let $b_{\bullet}^G = \dim H_{\bullet}(X;\R)^G$, where $\bullet = 1,2,+$.

Our main theorem is,
\begin{Theorem}\label{thm:main}
Let $G=\Z_p$ be the cyclic group of odd prime order $p$,
and $X$ a closed oriented $4$-manifold with $b_1=0$ and $b_+\geq 2$.
Suppose $G$ acts on $X$ with $b_+^G\geq 2$, and the $G$-action has a lift to a $\Spinc$-structure $c$ with $d(c)=0$.
If $d(c,G_\alpha)\leq 0$ for any lift of the $G$-action,
then
\begin{equation}\label{eq:modp-eq}
\SW(X,c)\equiv \sum_{\alpha}\m_\alpha\SW(X,c,G_\alpha)\mod p,
\end{equation}
where $\m_\alpha$ are integers determined by the $G$-index of the Dirac operator and the $G$-action on $H^+(X;\R)$. 
{\rom (}If $d(c,G_\alpha)< 0$, then $\m_\alpha =0$. 
For the other case, see around \eqref{eq:mult} for the precise definition.\rom{)}
\end{Theorem}
\begin{Remark}
\thmref{thm:main} can be generalized to the case when $p=2$ or $b_1>0$ with appropriate assumptions.
To avoid a complicated description, we only give the proof of the case of \thmref{thm:main}, and the detail of such generalizations will be left to readers.
Other possibilities of generalizations will be referred in \remref{rem:gen} below.
\end{Remark}

The strategy of the proof of \thmref{thm:main} is analogous to those in \cite{N} and \cite{Np}.
We will work out a $G$-equivariant perturbation of the monopole map. 
Under the $G$-action, the moduli space splits into two parts: the $G$-fixed part and the $G$-free part.
When the dimension of the moduli is $0$, the number of solutions in the $G$-free part is a multiple of $p$.
On the other hand, the number of $G$-invariant solutions {\it is} the $G$-monopole invariant.
However, the transversality is not necessarily achieved on these $G$-invariant solutions.
Then, we give a canonical way of $G$-equivariant perturbation, which enables us to determine the multiplicities of these solutions.

The organization of the paper is as follows:
Section 2 gives a brief review on $G$-monopole invariants.
In Section 3, we prove \thmref{thm:main}.
In Section 4, we discuss several examples.
%
%
\section{$G$-monopole invariants}\label{sec:G-monopole}
%
%
In this section, we give a brief review on $G$-monopole invariants defined by Ruan \cite{Ruan}.

Let $X$ be a closed oriented $4$-manifold, and $c$ a $\Spinc$-structure on $X$. 
Let $\G$ be the gauge transformation group which consists of automorphisms of $c$ covering the identity map of $X$.
Note $\G=\Map (X,S^1)$.
We introduce another automorphism group  $\tilde{\G}$ consisting of pairs $(f,\tilde{f})$, where $f\colon X\to X$ is a diffeomorphism of $X$, and $\tilde{f}\colon c\to c$ is an automorphism of $c$ covering $f$. 
Then, we have an exact sequence,
$$
1\to \G\to\tilde{\G}\to \Diff^+(X),
$$
where $\Diff^+(X)$ is the group of orientation-preserving diffeomorphisms of $X$.

Let $G$ be a finite group.
Note that giving an effective orientation-preserving $G$-action on $X$ is equivalent to giving a subgroup $G$ of $\Diff^+(X)$.
Suppose such a $G$-action on $X$ is given, and $c$ satisfies $g^*c\cong c$ for any $g\in G$. 
Then the following group extension exists,
\begin{equation}\label{eq:ext}
1\to\G\to\hat{\G}\overset{\theta}{\to} G\to 1.
\end{equation}
Note that giving a lift of the $G$-action to $c$ is equivalent to giving a slitting of \eqref{eq:ext}, that is, giving a subgroup $G_\alpha$ of $\hat{\G}$ which is isomorphic to $G$ via $\theta$.

Suppose we are given data $(X,c,G_\alpha)$ as above.
In such a situation, Y.~Ruan defined the {\it $G$-monopole invariant} \cite{Ruan} as follows. 
In this case, the Seiberg-Witten equations are $G_\alpha$-equivariant, and the $G_\alpha$-invariant moduli space $\M(X,c,G_\alpha)$ is defined as the set of equivalence classes of $G_\alpha$-invariant solutions modulo $G_\alpha$-invariant gauge transformations.
The virtual dimension of $\M(X,c,G_\alpha)$ is given by,
$$
d(c,G_\alpha) = 2\ind D^{G_\alpha} - (1-b_1^G+b_+^G),
$$
where $\ind D^{G_\alpha}$ is the virtual dimension of the trivial part of the $G_\alpha$-equivariant Dirac index.
Note that we can orient all of $\M(X,c,G_\alpha)$ at the same time by fixing an orientation of $(H^1(X;\R)\oplus H^+(X;\R))^G$.
If $d(c,G_\alpha) =0$, then the $G$-monopole invariant $\SW(X,c,G_\alpha)$ is defined as the signed count of the number of elements in $\M(X,c,G_\alpha)$.
In general, we need to perturb the equations to avoid reducibles and achieve transversality. 
The standard argument proves that $\SW(X,c,G_\alpha)$ is well-defined if $b_+^G\geq 2$, and depends on chambers if $b_+^G=1$.

%
%
\section{Proof of \thmref{thm:main}}\label{sec:proof}
%
%
In this section, we prove our main theorem(\thmref{thm:main}).

Suppose $(X,c)$ with a $G=\Z_p$-action satisfies the conditions in \thmref{thm:main}, and a lift of the $G$-action to $c$, say $G_0$, is given.
Fix a $G$-invariant metric and a $G_0$-invariant connection $A_0$ on the determinant line bundle of $c$.
Then the monopole map $\mu$ is a proper $G\times S^1$-equivariant map.
Taking a finite dimensional approximation of $\mu$ \cite{BF}, we have a $G\times S^1$-equivariant map between finite rank representations:
$$
f_0\colon V\oplus R\to W\oplus R\oplus H,
$$
where $V$ and $W$ are complex representations of $G$ on which $S^1$ acts by multiplication, and $R$ and $H=H^+(X;\R)$ are real representations of $G$ on which $S^1$ acts trivially.
More explicitly, when $\C_j$ is the complex $1$-dimensional weight $j$ representation of $G$, $V$ and $W$ can be written as,
\begin{align*}
V&=\C_0^{a_0}\oplus \C_1^{a_1}\oplus \cdots \C_{p-1}^{a_{p-1}},\\
W&=\C_0^{b_0}\oplus \C_1^{b_1}\oplus \cdots \C_{p-1}^{b_{p-1}}.
\end{align*}
The $G_0$-index of the Dirac operator associated to $A_0$ is written as
$$
\ind_{G_0}D_{A_0}=\sum_{j=0}^{p-1}(a_j-b_j)\C_j.
$$
Note that the other lifts $G_\alpha$ are parameterized by $\alpha=j$ where $1\leq j \leq p-1$, and each $G_j$ is obtained by twisting $G_0$ by multiplication of $e^{-2\pi\sqrt{-1}j/p}$.
In other words, as $G_j$-representations, $V$ and $W$ become $V\otimes \C_{-j}$ and $W\otimes \C_{-j}$.

First, perturb $f_0$ $G\times S^1$-equivariantly so that the zero locus does not contain any reducible as follows:
Take a nonzero element $v$ in $H^+(X;\R)^G$, and perturb $f_0$ to $f:=f_0+v$.
Then $(f^{-1}(0))^{S^1}=\emptyset$. 
(See \cite{Np}, Section 2.3.)

Dividing $f$ by $S^1$, we obtain a section $s\colon B\to E$ of the vector bundle $E\to B$ which is given by
\begin{gather*}
E=\left((V\setminus\{0\})\times R\right)\times_{S^1}(W\oplus R\oplus H),\\
B=(V\setminus\{0\})/S^1\times R.
\end{gather*}

When $d(c)=0$, $\SW(X,c)$ is the signed count of zero points of $s$ if $s$ is transversal to the zero section.
Note that $(V\setminus\{0\})/S^1$ is $G$-equivariantly homeomorphic to $P(V)\times \R_{+}$, where $P(V)$ is the projective space of $V$, and $\R_+$ is the space of positive real numbers.
The $G$-fixed point set of $P(V)$ can be written as (\cite{Np}, Lemma 3.1),
$$
P(V)^G = \coprod_{j=0}^{p-1}P(\C_j^{a_j}).
$$
Let $B_j = \R_+\times P(\C_j^{a_j})\times R_0$, where $R_0$ is the $G$-fixed part of $R$.
Then the $G$-fixed point set of $B$ decomposes into  its connected components as $B^G=B_0\cup B_1\cup\cdots\cup B_{p-1}$. 

Note that each $B_j$ corresponds to the lift $G_j$, and 
$$
d(c,G_j)=2(a_j-b_j)-(1+b_+^G).
$$
When $d(c,G_j)=0$, $\SW(X,c,G_j)$ is given by the signed count of zero points of $s|_{B_j}$ if $s|_{B_j}$ is transversal to the zero section in $E^G$.

Now, let us carry out the $G$-equivariant perturbation of $s$.
When $d(c,G_j)<0$, we can perturb $s$ $G$-equivariantly around $B_j$ so that $s^{-1}(0)\cap B_j=\emptyset$.
When $d(c,G_j)=0$, we can perturb $s$ $G$-equivariantly around $B_j$ so that $s|_{B_j}$ is {\it transversal along $B_j$}.
Then, the problem is how to count multiplicities of zero points on $B_j$. 

Let $x$ be a point in $s^{-1}(0)\cap B_j$.
We would like to describe the differentiation $(Ds)_x$ of $s$ at $x$.
The tangent space of $B$ at $x$ decomposes into the $G$-invariant direction and its complement:
$T_xB = T_xB_j\oplus V^\prime.$
Then $T_xB_j$ and  $V^\prime$ can be identified as
\begin{gather*}
T_xB_j= \R\times \C_j^{a_j}\times R_0,\\
V^\prime = \sum_{k\neq j}\C_{k-j}^{a_k}\oplus R^\prime,
\end{gather*}
where $R^\prime$ is the orthogonal complement of $R_0$ in $R$.
By reordering $\C_j$'s in $V^\prime$, rewrite  $V^\prime$ as
$V^\prime = \C_1^{a_1^\prime}\oplus\cdots\oplus\C_{p-1}^{a_{p-1}^\prime}\oplus R^\prime$,
where $a_k^\prime = a_{k+j}$.

Similarly, the vertical tangent space $V_{s(x)}E$ of $E$ at $s(x)$ decomposes as, $V_{s(x)}E = W_0\oplus W^\prime$, 
where $W_0$ is $G$-invariant part and $W^\prime$ is its complement.
When we decompose $H=H^+(X;\R)$ into $H_0\oplus H^\prime$, where $H_0$ is the $G$-fixed part and $H^\prime$ its complement, $W^\prime$ can be identified with
$$
W^\prime = \sum_{k\neq j} \C_{k-j}^{b_j}\oplus R^\prime\oplus H^\prime.
$$
Let us choose orientations of $H$ and $H_0$ (hence $H^\prime$ too), and fix an arbitrary identification $H^\prime = \C_1^{h_1}\oplus\cdots\oplus \C_{p-1}^{h_{p-1}}$ so that $H=H_0\oplus H^\prime$ and $H_0\oplus \C_1^{h_1}\oplus\cdots\oplus \C_{p-1}^{h_{p-1}}$ have same orientation.
(Here, we used the assumption that $p$ is odd.)
Rewrite $W^\prime$ as $W^\prime = \C_1^{b_1^\prime}\oplus\cdots\oplus\C_{p-1}^{b_{p-1}^\prime}\oplus R^\prime$,
where $b_k^\prime=b_{k-j} + h_k$.
Let $L_0$ be the linear map which is the composition of the following maps:
$$
\begin{CD}
L_0\colon V^\prime @>{Ds_x}>> T_{s(x)}E@>{p_v}>> V_{s(x)}E @>{p_w}>> W^\prime,
\end{CD}
$$
where $p_v$ and $p_w$ are the orthogonal projections.

We will {\it cancel out} common parts in $V^\prime$ and $W^\prime$ by a perturbation by a $G$-linear map.
We give a local model of this as follows.
Let $e_k =\min \{a_k^\prime,b_k^\prime\}$.
We can take an orientation-preserving  $G$-linear map $l\colon V^\prime\to W^\prime$ so that $\im(L_0+l)\cong\sum_k\C_k^{e_k}\oplus R^\prime$. 
Let $W_e=\im(L_0+l)$ and its complement in $W^\prime$ be $W_r$, 
and $V_r=\ker (L_0+l)$ and its complement in $V^\prime$ be $V_e$.
Then,
\begin{gather*}
V^\prime = V_e\oplus V_r,\\
W^\prime = W_e\oplus W_r,\\
V_e\cong W_e\cong \sum_k\C_k^{e_k}\oplus R^\prime.
\end{gather*}

Next, we give a local model of perturbation in the direction of $V_r$. 
Let 
$$
I = \{k\,|\,m_k = a_k^\prime-e_k>0\} \text{ and } I^\prime = \{k\,|\,n_k = b_k^\prime-e_k>0\}.
$$
Then 
$$
V_r=\sum_{k\in I}\C_k^{m_k},\ W_r =\sum_{k\in I^\prime}\C_k^{n_k}. 
$$
Note that $I\cap I^\prime=\emptyset$ and $\dim V_r = \dim W_r$.
We will perturb $s$ around $x$ by a (nonlinear) $G$-equivariant map $\psi\colon V_r\to W_r$.
The next example will illustrate how to take $\psi$.
\begin{Example}
Suppose $G=\Z_5$, $V_r=\C_1\oplus\C_4$ and  $W_r=\C_2\oplus\C_3$.
Then take $\psi\colon  \C_1\oplus\C_4\to \C_2\oplus\C_3$ which is given by 
$\psi(z,w)=(z^2,w^2)$. 
If we perturb $s$ around $x$ by $\psi$, then the multiplicity of $x$ is equal modulo $5$ to $2\times 2=4$. 
As another choice, we can take  $\psi$ given by $\psi(z,w)=(w^3,z^3)$.
In this case, the multiplicity of $x$ is also equal modulo $5$ to $3\times 3\equiv 4$. 
The multiplicity $4$ can be calculated by $2\cdot 3/1\cdot 4\equiv 4$ in the finite field $\F_5$.
\end{Example} 
The general case is given as follows.
Let $(z_1,\ldots, z_r)$ be the coordinate of $V_r$ where $z_k\in \C_{i_k}$,
and $(w_1,\ldots, w_r)$ be that of $W_r$ where $w_k\in \C_{i_k^\prime}$.
Then $\psi\colon V_r\to W_r$ is given by 
$$
\psi(z_1,\ldots,z_r)= (z_1^{i^\prime_1/i_1},\ldots, z_r^{i^\prime_r/i_r}),
$$
where $i^\prime_k/i_k$ is calculated in $\F_p$, and identified with an integer which represents it.

The multiplicity $\m_j$ of $x$ is given by
\begin{equation}\label{eq:mult}
\m_j = \frac{\prod_{k=1}^r i^\prime_k}{\prod_{k=1}^r i_k}.
\end{equation}

By using an appropriate $G$-invariant cut-off function, perturb the section $s$ around $x$ by $l + \psi$.
For every point in $s^{-1}(0)\cap B^G$, such a perturbation should be carried out.
We also need to perturb $s$ $G$-equivariantly on the free part $B\setminus B^G$.
This is easy.

Now, we complete the proof of \thmref{thm:main}.
By the perturbation so far, each of zeros of $s$ on $B^G$ has its multiplicity $\m_j$.
On the other hand, $G=\Z_p$ acts freely on $s^{-1}(0)\cap (B\setminus B^G)$.
Hence, the relation \eqref{eq:modp-eq} holds.

\begin{Remark}\label{rem:gen}
In the proof above, the assumption $d(c)=0$ is not essential.
In the case when $d(c)>0$, we can use the technique of {\it cutting down the moduli space} as in \S 3(iii) in \cite{Np}. 
On the other hand, the assumption $d(c,G_\alpha)\leq 0$ seems essential to our proof.
It  would be an interesting problem to consider the case when $d(c,G_\alpha)>0$. 
\end{Remark}
\begin{Remark}
Another possibility of generalization is to consider $p$-fold branched coverings.
As mentioned in \secref{sec:intro}, the case of $2$-fold branched covering is studied by \cite{RW,Park,Cho}.
One could try to prove similar results for higher orders. 
\end{Remark}

%
%
\section{Examples}\label{sec:examples}
%
%
In this section, we give several examples.
\subsection{Example 1}\label{subsec:ex1}
Let $X$ be the $K3$ surface of the Fermat type in $\CP^3$ defined by the equation $z_0^4+z_1^4+z_2^4+z_3^4=0$.
Let $G=\Z_3$ act on $X$ by permutation of components.
Let $c$ be the $\Spinc$-structure determined by the spin structure, and consider the lift $G_0$ of the $G$-action to $c$ whose induced action on the determinant line bundle is just the diagonal action $X\times \C_0$.
Then, the $G_0$-index of the Dirac operator is written as $\ind_{G_0} D = 2\C_0$. (See \cite{LN}.)
The finite dimensional approximation of the monopole map has the form,
$$
f\colon\C_0^{x+2}\oplus\C_1^y\oplus \C_2^z\to \C_0^{x}\oplus\C_1^y\oplus \C_2^z\oplus\R^3,
$$
where $\R$ is the real $1$-dimensional trivial representation.
It follows that $d(c,G_0)=0$, $d(c,G_1)=d(c,G_2)<0$, and, by \thmref{thm:main},
$$
\SW(X,c)\equiv \SW(X,c,G_0)\mod 3.
$$ 
In fact, $\SW(X,c)= \SW(X,c,G_0)=1$, because there exists the unique $G$-invariant solution with constant spinor by the perturbation by a $G$-invariant holomorphic $2$-form.
We remark that the action in Proposition 4.11 of \cite{LN2} gives a similar example in the case of $G=\Z_5$
%
%
\subsection{Example 2}\label{subsec:ex2}
D.-Q.~Zhang introduced a holomorphic $G=\Z_3$-action on a $K3$ surface $X$ with $b_+^G=1$ (\cite{Zh}, Example 5.3, due to S.~Tsunoda).
Let $c$ be the spin, and consider the lift $G_0$ as in \subsecref{subsec:ex1}.
In this case, the finite dimensional approximation is of the form, 
$$
f\colon\C_0^x\oplus\C_1^{y+1}\oplus \C_2^{z+1}\to \C_0^{x}\oplus\C_1^y\oplus \C_2^z\oplus\R\oplus \C_1.
$$
Then, $d(c,G_0)<0$, $d(c,G_1)=d(c,G_2)=0$.
Note that $b_+^G=1$ in this case, and therefore $\SW(X,c,G_\alpha)$ depend on chambers.
Nevertheless, the formula \eqref{eq:modp-eq} in \thmref{thm:main} holds for any chamber as
$$
\SW(X,c)\equiv \SW(X,c,G_1) + 2\SW(X,c,G_2)\mod 3.
$$ 

In fact, the following occurs:
\begin{Proposition}\label{prop:chamber}
In a chamber $\CC_+$, $\SW(X,c,G_1)=1$ and $\SW(X,c,G_2)=0$.
In another chamber $\CC_-$, $\SW(X,c,G_1)=0$ and $\SW(X,c,G_2)=-1$.
\end{Proposition}
\begin{Remark}
In the chamber $\CC_+$, the formula \eqref{eq:modp-eq} holds as $1\equiv 1+2\cdot 0$.
On the other hand, in $\CC_-$, the formula \eqref{eq:modp-eq} holds as $1\equiv 0+2\cdot (-1)$.
\end{Remark}
To prove \propref{prop:chamber}, we note the next.
\begin{Lemma}
$X$ admits a K\"{a}hler form $\omega$ preserved by the $G$-action.
\end{Lemma}
\proof
%
%
%
Let us recall the the construction of the log Enriques surface $\bar{S}=X/G$ (\cite{Zh}, Example 5.3).
Let $x$, $y$, $z$ be the homogeneous coordinates of $\CP^2$.
Consider three cuspidal cubic curves in $\CP^2$:
$$
C_1\colon x^3 =y^2z,\quad C_2\colon y^3 =z^2x,\quad C_3\colon z^3 =x^2y.
$$
Let $\xi$ be a primitive 7th root of the unity. 
Then $C_1\cap C_2\cap C_3 =\{(\xi^i:\xi:1)\,|\,0\leq i\leq 6\}$.
Let $\tau\colon S\to\CP^2$ be the blowing up of cusps $(1:0:0)$, $(0:1:0)$, $(0:0:1)$, and $7$ points in $C_1\cap C_2\cap C_3$.
Then $S$ contains three disjoint nonsingular $(-3)$-curves from $C_1$, $C_2$ and $C_3$.
Collapsing these $(-3)$-curves, we obtain the surface $\bar{S}$ whose covering is a $K3$.
These surfaces fit into the following diagram:
$$
\begin{CD}
X\#3\overline{\CP}^2@>{\pi}>>S@>{\tau}>>\CP^2\\
@V{\sigma}VV @V{c}VV\\
X@>{\bar{\pi}}>>\bar{S},
\end{CD}
$$
where $\sigma$ and $\tau$ are blowing up, $c$ is the collapsing map, $\pi$ is a $G$-fold covering branched along the $(-3)$-spheres, and $\bar{\pi}$ is a $G$-cover.
Note that $X\#3\overline{\CP}^2$ has a $G$-invariant K\"{a}hler form obtained by pulling back a K\"{a}hler form on $\CP^2$ via $\tau$ and $\pi$. 
By blowing down, we have a K\"{a}hler form $\omega$ on $X$ which is preserved by the $G$-action.
\endproof

\proof[Proof of \propref{prop:chamber}]
The positive spinor bundle $S^+$ of $c$ can be written as $S^+=I\oplus K_X^{-1}$, where $I$ is a trivial bundle and $K_X$ is the canonical line bundle of $X$ (which is also trivial).
Therefore, a spinor $\phi$ has two components $\phi=(\alpha,\beta)$. 
Since the $G$-action on $K_X^{-1}$ is given by $X\times \C_2$ and we fix the lift $G_0$ so that $\det S^+=I\otimes K_X^{-1} = X\times \C_0$, the $G_0$-action on $I$ is given by $I=X\times\C_1$.
By Taubes' perturbation \cite{Ko}({\it cf.} \cite{Taubes}) adding $-ir\omega$, we have a unique solution  such that $\alpha=\text{const.}$ and $\beta=0$. 
This solution is $G_1$-invariant.
On the other hand, if we use the perturbation adding $+ir\omega$, then the roles of $\alpha$ and $\beta$ are exchanged.
Therefore, we have a unique solution with $\alpha=0$ and $\beta=\text{const.}$ 
which is $G_2$-invariant.
These two belong to different chambers.
By considering the orientations readily, the proof is completed.
\endproof
\begin{Remark}
In this case, the formula  $\SW(X,c)= \SW(X,c,G_1) - \SW(X,c,G_2)$ holds.
In fact, the perturbation adding $+ir\omega$ corresponds to a linear but orientation-reversing  perturbation by $\psi\colon \C_1\to\C_2$ given by $\psi(z)=\bar{z}$.
\end{Remark}

Several actions of higher order $G$ in \cite{Zh} give similar examples.


\begin{thebibliography}{99}





\bibitem{BF} S.~Bauer and M.~Furuta, 
{\it A stable cohomotopy refinement of Seiberg-Witten invariants: \Romnum{1}},
Invent. Math. {\bf 155} (2004), 1--19.


\bibitem{Cho} Y.~S.~Cho and Y.~H.~ Hong,
{\it Seiberg-Witten invariants and {\rom (}anti-{\rom )}symplectic involutions},
Glasgow Math.~J. {\bf 45} (2003), 401--413.

\bibitem{Fang} F.~Fang,
{\it Smooth group actions on $4$-manifolds and Seiberg-Witten invariants},
 Internat. J. Math. {\bf 9}, No.8 (1998) 957--973. 

\bibitem{Ko} D.~Kotschick,
{\it The Seiberg-Witten invariants of symplectic four-manifolds {\rom (}after C. H. Taubes{\rom)}},
S\'{e}minaire Bourbaki, Vol. 1995/96.
Ast\'{e}risque {\bf 241} (1997), Exp. No. 812, 4, 195--220. 

\bibitem{LN} X.~Liu and N.~Nakamura, 
{\it Pseudofree $\Z/3$-actions on $K3$ surfaces}, 
Proc. Amer. Math. Soc. {\bf 135} (2007), no. 3, 903--910. 

\bibitem{LN2} X.~Liu and N.~Nakamura, 
{\it Nonsmoothable group actions on elliptic surfaces}, 
Topology Appl., {\bf 155} (2008), no.9, pp 946-964.


\bibitem{N} N.~Nakamura,
{\it A free $\Zp$-action and the Seiberg-Witten invariants},
J.~Korean Math.~Soc. {\bf 39} (2002), No. 1, 103--117.

\bibitem{Np} N.~Nakamura,
{\it Mod $p$ vanishing theorem of Seiberg-Witten invariants for $4$-manifolds with $\Z_p$-actions},
Asian J. Math. {\bf 9}, (2006), no.4, 731--748.

\bibitem{Park} B.~D.~Park,
{\it Seiberg-Witten invariants and branched covers along tori},
 Proc. Amer. Math. Soc.  {\bf 133}  (2005),  no. 9, 2795--2803. 

\bibitem{Ruan} Y.~Ruan,
{\it Virtual neighborhoods and the monopole equations,}  
Topics in symplectic $4$-manifolds (Irvine, CA, 1996),  101--116,
First Int. Press Lect. Ser., I, Int. Press, Cambridge, MA, 1998. 

\bibitem{RW} Y.~Ruan and S.~Wang,
{\it Seiberg-Witten invariants and double covers of $4$-manifolds},
Comm. Anal.~Geom. {\bf 8} (2000), No.3, 477-515.

\bibitem{Sz} M.~Szymik,
{\it Bauer-Furuta invariants and Galois symmetries},
preprint.

\bibitem{Taubes} C.~H.~Taubes,
{\it ${\rm SW}\Rightarrow{\rm Gr}$: from the Seiberg-Witten equations to pseudo-holomorphic curves},
J. Amer. Math. Soc. {\bf 9} (1996), no. 3, 845--918.  

\bibitem{Zh} D.-Q.~Zhang,
{\it Logarithmic Enriques surface}.
J.~Math.~Kyoto Univ. {\bf 31} (1991), 419--466.

\end{thebibliography}
\end{document}